\documentclass[twoside]{article}

\usepackage{sfuEng}

\usepackage{latexsym,amscd,amssymb,amsopn,amsthm,amsmath,enumerate}
\usepackage{cite}
\usepackage{} 

\makeatletter
\long\def\@makecaption#1#2{%
  \vskip\abovecaptionskip
  \sbox\@tempboxa{#1. #2}%
  \ifdim \wd\@tempboxa >\hsize
    #1. #2\par
  \else
    \global \@minipagefalse
    \hb@xt@\hsize{\hfil\box\@tempboxa\hfil}%
  \fi
  \vskip\belowcaptionskip}
\makeatother

\god{2011}

\nomer{4(3)}

\pp{320--331}                 

\firstvar{18.01.2011}         

\lastvar{25.02.2011}          

\toprint{10.03.2011}          

\udk{517.55}

\tit{Solving two dual problems of splicing vortex and potential flows with Goldshtik's variational method}

\abstract{The general problem of a perfect incompressible fluid motion with vortex areas and variant constant vorticities is formulated.
The M.A.~Goldshtik's variational approach is considered on research of dual problems for flows with vortex and potential areas that describe detached flow and a motion model of a perfect incompressible fluid in field of Coriolis forces.}

\keywords{vortex and potential flows, variational method,Green's function, extremum of the functional.}

\begin{document}

\maketitBegin 

\def\author{
\begin{flushright}
\large \bfseries  Isaac I.\, Vainshtein\footnote{\noindent\hbox{\parbox[t]{10cm}{\mbox{isvain@mail.ru}}}}

{\small \mdseries Institute of Space and Information Technologies,\\
Siberian Federal University,\\
Kirensky st., 26, Krasnoyarsk, 660041,\\
Russia }
\end{flushright}
}
\author{}

 \maketitEnd  


\markright{\footnotesize
    Isaac I.\,Vainshtein
    \hfill Solution of two dual problems of splicing vortex and potential flows ... }
%
%

\section{General formulation of a perfect incompressible fluid problem with vortex areas and constant vorticities}

A lot of publications and monographs cover researching of vortex flows of perfect incompressible fluid.
This subject is represented in every hydrodynamic course.
The M.A.~Goldshtik's monograph <<Vortex flows>> [1] is fundamental for current research.
Goldshtik's publication investigates vortex flows
that have areas which bounded with the same <<zero>> streamline with circled motion inside,
whereas external flow could be potential.
In the same monograph different examples of such flows in nature and engineering are mentioned and
a research of problems that could be profitable in pure and applied science is represented.
These are detached flows and flows behind high-drag objects.
A field of centrifugal accelerations is an engineering appliance of rotating flows.
As considering machines could be met various separators, dryers, combustion chambers,
cyclonic separations (dust filters),  vortex nuclear reactors and etc.

Steady-state vortex flow of a perfect incompressible fluid in plain case could be described by equation below.
$$
\Delta\Psi=F(\Psi)
$$
where vorticity $\omega=F(\Psi)$ is an arbitrary composite function of stream function $\Psi$,
$v_x=\dfrac{\partial\Psi}{\partial y}$,
$v_y=-\dfrac{\partial\Psi}{\partial x}$. Vorticity it self satisfies equation
$$
 \frac{d\omega}{dt}=\frac{\partial\omega}{\partial t}+v_x \frac{\partial\omega}{\partial x}+
 v_y \frac{\partial\omega}{\partial y}= \nu\Delta\omega,
$$
where $\nu$ is a kinematic viscosity [1]. Besides, in steady-state case equation for vorticity doesn't contain $\omega$ value.
The bilateral maximum principle for $\omega$ value is valid in this case.
Mentioned above allows to prove the assumption of vortex constancy (given case is considered in this research) for
limiting steady-state  motion of viscous fluid in domain that bounded with closed streamline when viscosity $\nu$ vanishes [1--3].

 This statement leads to a consideration of flow categories where areas contain closed streamlines of various steady-state vorticities.
 Even it is possible to suggest that number of described areas is unlimited.
 Volumes of vorticities could be either arbitrary sign or vanished.
 In the last case the flow in area is potential one.

Thus, we are going to a problem statement of a perfect incompressible flow with piecewise-constant vorticities that generalizes problems (1), (3) and (2), (3) considered below.
Vorticities $\omega_1,\omega_2, \ldots, \omega_n$ are given
and $D$ is a domain of the stream with bound $\Gamma$ (domain $D$ could be unbounded).
It's required to construct stream areas $B_i, \
(i=1,2,\ldots,n)$, $\bigcup B_i=D$ and find
continuously differentiable in domain $D$ stream function $\Psi(x,y)$ that in every domain
$B_i$ satisfies the equation $\Delta \Psi =\omega_i$. The function $\Psi(x,y)=\varphi(s)$ at the border $\Gamma$ of stream domain $D$ whereas function $\Psi(x,y)$ is equal to the same constant value at the borders of each area $B_i$.

In order to do a complete research of possible flows it's needed to include
to the problem statement a possibility of arbitrarily vanishing vorticity values to consider the existence of areas with potential flows.

Considered in this research equations (1), (2) specify flows with two areas in one of which a vorticity is nonzero whereas in the other it vanishes.
Stream function vanishes at the areas' bounds.

Mentioned above leads to two dual problems of splicing vortex and potential
flows [1,\,4--6]. Conditions for stream function's sign limit the possible flow categories.

As an example of the flow with potential and two vortex areas could be a flow around rectangular trench.
The flow could divide into three areas as it's stated in monograph of M.A.~Lavrentyev, B.V.~Shabat [7].

A motion is potential in the first area which is located over the trench.
A motion has steady-state vorticities in two other areas which are located in the trench one over another:
$-\omega$ in upper area and $\omega$ in lower one. In the publication [8] there is considered a given problem case for two vortex areas with $\omega_1$ and $-\omega_2$ vorticities where $\omega_1$ and
$\omega_2$ are nonzero constants:
$$
\Delta\Psi=\left\{ \begin{array}{ll} \displaystyle
\omega_1,\quad\mbox{if}\ \Psi<0, \\[1mm]
\displaystyle -\omega_2,\quad\mbox{if}\ \Psi>0,
\end{array}\right.
$$
$$
\Psi|_\Gamma=\varphi(s)\geq 0,
$$
and existence conditions of these flows are derived.

\section{Dual problems of splicing vortex and potential flows. \\ M.A.~Goldshtik's variational approach}

Two dual problems are known in theory of vortex and potential flows:
in bounded plain domain $D$ with bound $\Gamma$
it is required to find continuously differentiable solutions of equations
\begin{equation}\label{1}
\Delta\Psi=\left\{ \begin{array}{ll} \displaystyle
\omega,\quad\mbox{if}\ \Psi<0, \\
\displaystyle  0,\quad\mbox{if}\ \Psi>0,
\end{array}\right.
\end{equation}
\begin{equation}\label{2}
\Delta\Psi=\left\{ \begin{array}{ll} \displaystyle
\omega,\quad\mbox{if}\ \Psi>0, \\
\displaystyle  0,\quad\mbox{if}\ \Psi\leq0
\end{array}\right.
\end{equation}
with boundary condition
\begin{equation}\label{3}
\Psi|_\Gamma=\varphi(s)\geq0,\quad \omega \mbox
{(vorticity)}>0
\end{equation}
where $\Psi(x,y)$ is a stream function [1, 6, 7, 8]. In what follows it is supposed that
bound $\Gamma$ is piecewise and function
$\varphi(s)$ is that there exists a harmonic in domain $D$ function that possesses the value $\varphi(s)$ at a bound.

The problem (1), (3) describes M.A.~Lavrentyev's scheme of detached flows according to which
a detach area is separated from external flow with zero streamline.
It is considered that flow is vortex inside the area with constant vorticity $\omega$.
The flow is potential outside the area [1,\,4].

The problem (2), (3) describes a model of flat motion of perfect fluid in field of Coriolis forces  [1,\,5,\,6].

Function $\Psi_0(x,y)$ denotes a harmonic function that satisfies condition (3).
This function is a solution of problem (1), (3) and positive in domain $D$.
Described above solution is named trivial.
In [4] proved the existence of nontrivial (with negative area)
solution of problem (1),(3) for reasonably great value $\omega$. In [6,\,8] derived the following inequality
\begin{equation}\label{4}
\omega>\frac{4Ce}{R^2},\quad C=\max\varphi(s),
\end{equation}
when nontrivial solution exists. Here $R$ is a radius of the area-largest inscribed in domain $D$ circle.

\begin{figure}[h!]
\centering
\includegraphics [width=4.2cm] {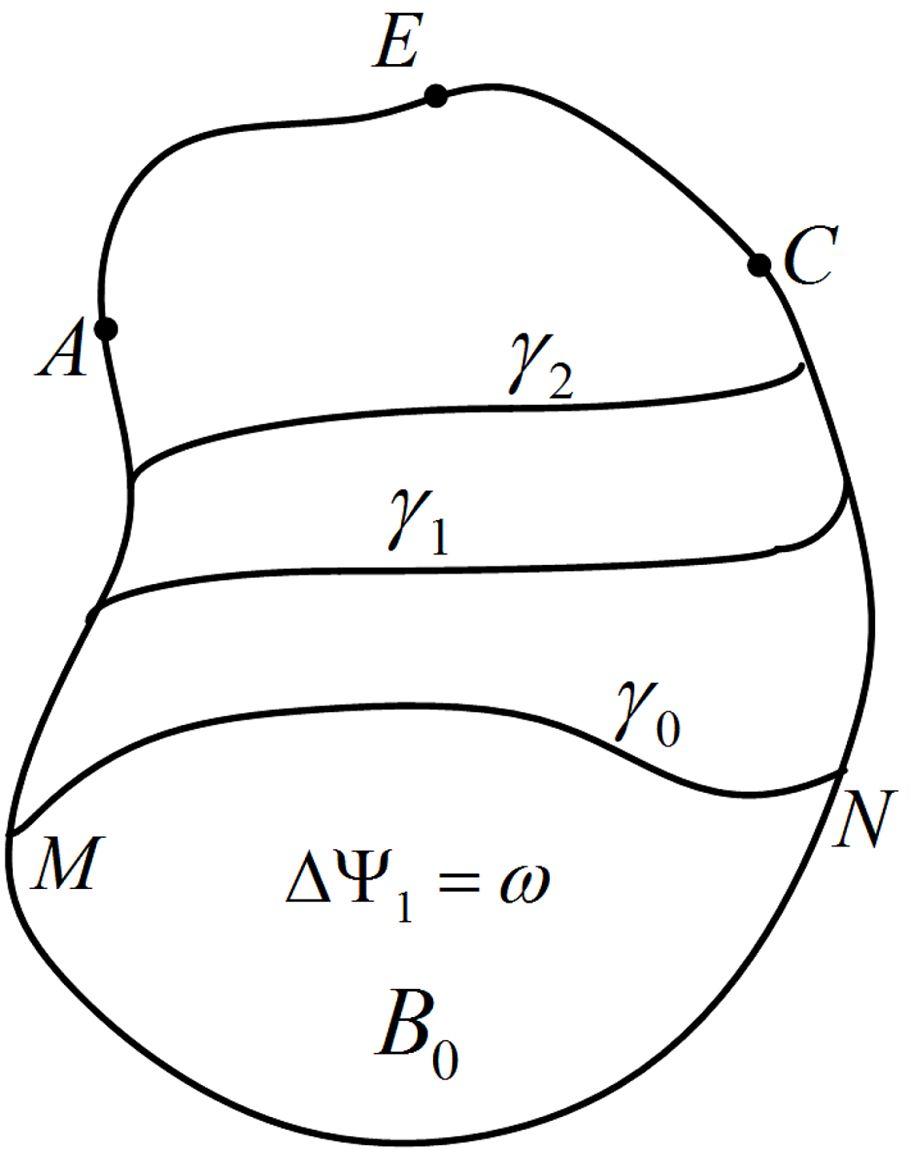}
\caption{Flowing area scheme}
\end{figure}

Let's introduce the scheme of M.A.~Goldshtik's proof [1,\,4] of nontrivial solution existence.
Let $\varphi(s)$ be nonzero only at
$AEC$ piece of bound $\Gamma$ (Fig.1). The case when $\varphi(s)>0$ along the hole bound
 $\Gamma$ is included. A smooth curve $\gamma_0$ is crossing $AEC$ and bounding with piece of $\Gamma$ domain
$B_0$.
The following function is constructed.
\begin{equation}\label{5}
\Psi_1(x,y)= \Psi_0(x,y)-\frac{\omega}{2\pi}\int\limits_{B_0}\!\!\!\int G(x,y,\xi,\eta)d\xi d\eta,
\end{equation}
where $G(x,y,\xi,\eta)$ is a Green's function for Laplace operator of Dirichlet problem in domain $D$.
The volume $\omega$ satisfies inequality
\begin{equation}\label{6}
\omega>\frac{\Psi_0}{\frac{1}{2\pi}\int\limits_{B_0}\!\!\!\int G(x,y,\xi,\eta)d\xi d\eta}
\end{equation}
for all points of curve $\gamma_0.$ So function $\Psi_1(x,y)$
is negative at $\gamma_0$ and because of boundary condition $\Psi_1>0$ nearby $AEC$, there exists curve $\gamma_1$ (contour curve of harmonic function) at which $\Psi_1=0.$
The curve $\gamma_1$ lies out of $B_0$ and rests with it's end to curves  $AM$ and $CN$ of bound $\Gamma.$
Then we define function
$$
\Psi_2(x,y)= \Psi_0(x,y)-\frac{\omega}{2\pi}\int\limits_{B_1}\!\!\!\int G(x,y,\xi,\eta)d\xi d\eta,
$$
where $B_1$ is a domain of negative function $\Psi_1$, $B_1\supset
B_0$. At the curve $\gamma_1$ (except for ends) takes place the following inequality $\Psi_2<0.$
Similarly to previous statement the curve $\gamma_2$ and domain $B_2\supset
B_1$ is defined. Continuing described process  we construct a sequence of functions
\begin{equation}\label{7}
\Psi_n(x,y)= \Psi_0(x,y)-\frac{\omega}{2\pi}\int\limits_{B_{n-1}}\!\!\!\!\!\int G(x,y,\xi,\eta)d\xi d\eta
\end{equation}
and domains $B_n$, $B_n\supset B_{n-1}.$ Now we proof that a consequence of functions $\Psi_n$ convergents to
solution of problem (1), (3).

If write equation (1) in form
$$
\Delta \Psi=\frac{\omega}{2}(1-\textrm{sign}\Psi),
$$
then offered in [1,\,4] method of finding nontrivial solution could be written as an iteration scheme

\begin{equation}\label{8}
\Delta \Psi_n=\omega(1-\textrm{sign}\Psi_{n-1}),\quad \Psi_n|_\Gamma=\varphi(s),
\end{equation}
where $\Psi_1$ is defined by formula \eqref{5}.

To model a perfect fluid under field of Coriolis forces for problem (2), (3)
it is required the fulfilment of inequality [1,\,5]
\begin{equation}\label{9}
\omega>\min\limits_{D}\frac{2\pi\psi_0}{\int\limits_D \!\!\!\int G(x,y,\xi,\eta)d\xi d\eta}=L.
\end{equation}

The domain $D$ is divided into two areas: circulating where $\Psi>0$ and
uncirculating where fluid is still ($\Psi\equiv 0).$ In the case of
$\omega<L$ an uncirculating area is absent and it could be considered as trivial solution.
The existence and uniqueness of solution of (2), (3) is proved in publications [1,\,5].

Scheme of the problem solution proof consists of a problem approximation with nonlinear discontinuity by
sequence of nonlinear problems with continued nonlinearity in right side of equation.

$$
\Delta \Psi_n=\frac{\omega}{2}\tanh(\Psi_nn),\quad \Psi_n|_\Gamma=\varphi(s).
$$

Existence of solution $\Psi_n$ for each $n$ and
convergence of sequence $\Psi_n$ to solution of problem (2), (3) is proving.

In the model problem when domain $D$ is a circle with radius $R$ and
$\varphi(s)=C$ when satisfying \eqref{4}  the problem (1),(3)
has two nontrivial solutions, there is one solution in the equality case and there are no nontrivial solutions for $\omega<4Ce/R^2$ [1]. The solutions could be written in explicit form. The existence of two nontrivial solutions is proved in [1,\,9] if domain $D$ is closely approximated to circle and boundary function is close to constant.

The differences of solutions characteristics of dual problems
could be explained by different behavior of right sides of equations (1), (2).
The function in right side of equation (1) is decreasing whereas in equation (2) is not decreasing by $\Psi.$

The existence proof of nontrivial solution of problem (1), (3)
could be found in publications [10--13]. Problems (1), (3) and (2), (3) for
arbitrary vorticity $\omega=F(\Psi)$ are considered in publications [14--16]
where for model examples with vorticity
$\omega=e^{\lambda\Psi}$ was established the increasing number of solutions effect in comparison with constant vorticity problems. Problem of splicing vortex and potential
flows in unbounded domain was considered in publications [17--19].
General case of vorticity $\omega=F(\Psi)$ was considered in [19].

Analog of problem (1), (3) of splicing vortex and potential flows in axisymmetric case
was considered in publications [6,\,20].

It should be mentioned that problems of number of nontrivial solutions and
area-dividing line properties for problem (1), (3) in general case are still unsolved.
The importance of this problems solution is highlighted in monographs [1,\,7].

To study the question about number of solutions in general case of problem (1), (2)
M.A.~Goldshtik suggested a variational formulation of a problem .
The following functional is considered.
\begin{equation}\label{10}
I(B)=\int\limits_{D}\!\!\!\int(\nabla \Psi)^2dxdy +
2\omega\int\limits_{B}\!\!\!\int\Psi dxdy,
\end{equation}
\begin{equation}\label{11}
\Psi(x,y)=\Psi_0(x,y)-\frac{\omega}{2\pi}\int\limits_{B}\!\!\!\int G(x,y,\xi,\eta)d\xi d\eta.
\end{equation}
The functional is defined on the aggregate of closed sets $B\subset D$ and could be transformed to
\begin{equation}\label{12}
I(B)=\int\limits_{\Gamma}\Psi_0\frac{\partial\Psi_0}{\partial n}ds +
\omega\int\limits_{B}\!\!\!\int(\Psi+\Psi_0)dxdy.
\end{equation}
Existence of the first integral is required in (12).
The formula of functional increment is derived
\begin{equation}\label{13}
\Delta I=2\omega\int\limits_{\Delta}\!\!\!\int\Psi dxdy -
(2\pi)^{-1}\omega^2\int\limits_{\Delta}\!\!\!\int dxdy
\int\limits_{\Delta}\!\!\!\int G(x,y,\xi,\eta)d\xi d\eta,
\end{equation}
where $\Delta$ is an increment of domain $B$.

Using formula \eqref{13} it could be proved that if extremum of the functional is realized in set $B$
function $\Psi$ will vanish at bound $\gamma$ of the set $B$.

Then $I(B)$
is continuous and upper and lower limited.
It is proved that nontrivial solution of a problem
which we get in [1] with iteration method  \eqref{8},
realizes functional minimum. The domain $B$ where solution is negative is derived by extending domains $B_n$
where $\Psi_n\leq 0$ at bounds and positive variations decrease the functional.
Notice that if in a certain set $I(B)$ maximum is realized then the set doesn't contain inner points,
because elimination of sufficiently small neighborhood of a point $(x_0,y_0)$ in which $\Psi<0$ increase the functional.
M.A.~Goldshtik stopped at this stage and advised to find conventionally critical points of the functional n order to prove second nontrivial solution.

The main purpose of this research is to spread variation approach that was suggested by M.A.~Goldshtik for problem (1), (3)
to dual problems (1), (3) and (2), (3).

\section{Solving of model dual problems with variational method}

Firstly we consider model dual problems (1), (3) and (2), (3)
when domain $D$ is a circle with radius $R$ and center in origin of coordinates
and boundary function $\varphi(s)=C.$

Solutions only depended on $r$ are searched.

Let's consider problem (1), (3) and the functional \eqref{12} for functions
\begin{equation}
\Psi(x,y)=C-\frac{\omega}{2\pi}\int\limits_{B_a}\!\!\!\!\int G(x,y,\xi,\eta)d\xi d\eta,
\end{equation}
where $B_a$ is a circle with radius $a$ and a center in origin of coordinates, $0\leq
a\leq R$. Then
\begin{equation}\label{15}
\Psi(x,y)=\left\{ \begin{array}{ll} \displaystyle \frac{\omega}{4}
(r^{2}-a^{2})+\frac{\omega a^{2}}{2}\ln\frac{a}{R}+C,
\quad\mbox{if}\  r\leq a, \\[2mm]
\displaystyle  \frac{\omega
a^{2}}{2}\ln\frac{r}{R}+C,\quad\mbox{if}\ a\leq r\leq R.
\end{array}\right.
\end{equation}
The function \eqref{15} is continuously differentiable. Satisfying condition
$\Psi=0$ when $r=a$ we get equation that define volume $a$:
\begin{equation}\label{16}
\frac{\omega a^{2}}{2}\ln\frac{a}{R}+C=0.
\end{equation}
Let's input function
\begin{equation}\label{17}
y(a)=\frac{\omega a^{2}}{2}\ln\frac{a}{R}+C,
\end{equation}
redefining it in the zero by limited value $C$.
Then
$$
y(0)=y(R)=C>0,\quad y^{\prime}(a)=\omega a\ln\frac{a}{R}+\frac{\omega a}{2}.
$$
A derivative of the function $y(a)$ in interval $(0,R)$ vanishes only in one point $a^*=R/\sqrt{2}.$
If conditions
$y(a^*)=-\dfrac{\omega R^2}{4e}+C<0 \ \left(\omega>\dfrac{4Ce}{R^2}\right)$
are satisfied and considering $y(0)=y(R)=C>0$ then equation \eqref{16} has two solutions
$a_1,\ a_2$,\quad $0<a_1<a_2<R.$ If increase $\omega$ a volume
$a_1$ tends to zero whereas $a_2$ tends to $R$. Thus, if inequality
$$
\omega>\frac{4Ce}{R^2}
$$
is satisfied then a model problem (1), (3) has two nontrivial solutions, but in the case of equality it has one solution.


Let's recur to functional \eqref{12}. After the direct integration of ($\Psi_0=C$) we will derive
\begin{equation}\label{18}
I(a)=\int \limits_{D}\Psi_0\frac{\partial \Psi_0}{\partial n}ds
+\omega\int\limits_{B_a}(\Psi+\Psi_0)dz =2\pi \omega a^2\left(\frac{\omega
a^2}{4}\ln\frac{a}{R}-\frac{\omega a^2}{16}+C\right).
\end{equation}
Let's research the functional $I(a).$ Now we have $I(0)=0,\quad
I(R)=2\pi \omega R^2\left(C-\dfrac{\omega R^2}{16}\right).$ From \eqref{18} it is followed
that for every fixed $\omega$ there is neighborhood $a>0$ where
$I(a)>0.$ Let's vanish the first derivative:
\begin{equation}\label{19}
I^{\prime}(a)=2\pi \omega a\left(\frac{\omega
a^2}{2}\ln\frac{a}{R}+C\right)=2\pi \omega ay(a)=0.
\end{equation}
Derived above equation when $a\neq 0$ coincides with equation \eqref{16}
of searching volume $a$ when solving the problem itself.
If inequality \eqref{4} is true then equation \eqref{19} has
two solutions $a_1,\ a_2$ $(a_1<a_2)$. There is a trivial solution in case $a=0$.

As in some neighborhood $a>0$ the functional $I(a)>0$ then in point $a_1$
the functional has absolutely positive maximum $I(a_1)>0$ and there is a minimum in point $a_2$.
Let the following inequality be true
\begin{equation}\label{20}
\omega>\frac{16C}{R^2},
\end{equation}
then $I(R)<0$ and condition \eqref{4} is fulfilled.
There is an absolute maximum in point $a_1$ whereas there is an absolute minimum $I(a_2)<I(R)<0$ in point $a_2$.

Thus, if we a searching the extremum of the functional in set of circles
$r\leq a\leq R$ and  inequality \eqref{4} is fulfilled then the functional has two extremums:
absolute maximum when $a=a_1$ and
absolute minimum when $a=a_1$. When substituting $a_1$ and
$a_2$ the function \eqref{15} gives two nontrivial solutions of model problem (1), (3).
In the case of equality in \eqref{4} there is one extremum and one nontrivial problem solution.

Let's consider a problem (2), (3) and a functional \eqref{12} for functions
\begin{equation}
\Psi(x,y)=C-\frac{\omega}{2\pi}\int\limits_{B_a}\!\!\!\int G(x,y,\xi,\eta)d\xi d\eta,
\end{equation}
where $B_a$ is a ring $a\leq r\leq R$. In this case
\begin{equation}\label{22}
\Psi(x,y)=\left\{ \begin{array}{ll} \displaystyle \frac{\omega}{4}
(a^{2}-R^{2})-\frac{\omega a^{2}}{2}\ln\frac{a}{R}+C,
\quad\mbox{if}\  r\leq a, \\[2mm]
\displaystyle  \frac{\omega}{4}
(r^{2}-R^{2})-\frac{\omega a^{2}}{2}\ln\frac{r}{R}+C,\quad\mbox{if}\ a\leq r\leq R.
\end{array}\right.
\end{equation}
The function \eqref{22} is continuously differentiable. Satisfying condition
$\Psi=0$ when $r=a$ an equation of finding $a$ is derived:
\begin{equation}\label{23}
\frac{\omega}{4}
(a^{2}-R^{2})-\frac{\omega a^{2}}{2}\ln\frac{a}{R}+C=0.
\end{equation}
Let's research equation \eqref{23}. We input function
\begin{equation}\label{24}
y_1(a)=\frac{\omega}{4}
(a^{2}-R^{2})-\frac{\omega a^{2}}{2}\ln\frac{a}{R}+C,
\end{equation}
and redefine in zero point with it's limit value. Now we have
\begin{equation}\label{25}
y_1(0)=C-\frac{\omega R^2}{4},\quad y_1(R)=C>0,\quad y_1^{\prime}(a)=-\omega a\ln\frac{a}{R}>0.
\end{equation}
The function $y_1(a)$ increases and, as followed from \eqref{24}, when
the inequality below is true
\begin{equation}\label{26}
\omega>\frac{4C}{R^2}
\end{equation}
it has only one root $a^*$ in interval $[0,R]$ (function
$y_1(a)$ has different signs in the ends of $[0,R]$).

Thus, the problem (2), (3) in the circle case and $\varphi (s)=C$ when inequality \eqref{26} is fulfilled, has only one solution \eqref{22} where $a=a^*$ is a root of equation \eqref{23}.

Let's recur to the functional \eqref{12}. After direct integration ($\Psi_0=C$) we get
\begin{equation}\label{27}
I(a)=2\pi\omega(C(R^2-a^2)+\frac{\omega a^4}{4}\ln\frac{a}{R}
+\frac{\omega R^2a^2}{4}-\frac{\omega R^4}{16}-\frac{3\omega a^4}{16}).
\end{equation}
Now let's research the functional $I(a).$ We have $I(0)=2\pi \omega
R^2(C-\dfrac{\omega R^2}{16}),\quad I(R)=0.$ The first derivative of the functional is vanished:
\begin{equation}\label{28}
I^{\prime}(a)=-4\pi\omega a(\frac{\omega}{4}
(a^{2}-R^{2})-\frac{\omega a^{2}}{2}\ln\frac{a}{R}+C)=-4\pi\omega ay_1(a)=0.
\end{equation}
Derived equation when $a\neq 0$ coincides with equation \eqref{23}
for searching volume $a$ when solving the problem itself.
When inequality \eqref{26} is fulfilled, there is one solution $a^*$.
When transiting through point $a^*$ a function $y_1(a)$ changes sign from minus to plus.
The functional's derivative changes sign from plus to minus according to \eqref{28}.
There is a maximum in point $a^*$. There are no other points of extremum when
satisfying inequality \eqref{26} $I(0)<0$. Whereas there is only one extremum
then $I(a^*)>0.$ The functional has absolute maximum.

Thus, if we search extremum of the functional in the set of rings that side to boundary of a circle
then satisfying inequality \eqref{26} there is only extremum (absolute maximum) for $a=a^*$ in the functional.
When $a=a^*$ the function \eqref{22} is a solution of model problem (2), (3) and realizes functional's extremum.

\section{The research of dual problems in general case}

Let's recur to general case of variational statement of the problem  (1), (2).
Assume that extremum of the functional is reached in domain $D$.
Now we show that there couldn't be a minimum in this case. Take at bound $\Gamma$ a point $M_0$ where $\Psi>0$.
Let $K(M_0,\rho)$ be a circle with radius $\rho$ and center at point $M_0$. The radius $\rho$ is chosen in order that in
$\Delta=D\cap K(M_0,\rho)$ the function $\Psi$ in form \eqref{11} is
greater than zero ($mes \Delta\neq 0).$ Let's consider increment of the functional
$\Delta I=I(D\backslash\Delta)-I(D).$ Taking into account that in $\Delta$
function $\Psi$ is positive we choose a sufficiently small $\rho$ (together with $mes
\Delta$) that according to \eqref{13} $\Delta
I$ will be negative. Described above means that exclusion of $\Delta$ from domain $D$ decreases the functional.

Now we consider a possibly of maximum. Let $\Psi>0$ in all points of domain $D$.
Function $\Psi$ couldn't be a solution of problem (1), (3) in considered case. Let's assume that $\omega$ satisfies inequality~\eqref{9}.
Then in domain $D$ there are points in which function $\Psi$ is negative.
Therefore $\Psi$  couldn't be a solution of problems (1),
(3) and (2), (3).

Let's estimate the volume $L.$ Let $R$ be a radius of the area-largest inscribed in
domain D circle (it could be considered that it's center is in origin of coordinates) and $G_R$ is a Green’s function for Laplace operator of Dirichlet problem in circle $C_R$. Let's consider function
$$
W=\frac{\omega}{2\pi}\int\limits_{D}\!\!\!\!\int
G(x,y,\xi,\eta)~d\xi
~d\eta-\frac{\omega}{2\pi}\int\limits_{C_R}\!\!\!\!\int
G_R(x,y,\xi,\eta)~d\xi ~d\eta.
$$
Described function is harmonic inside circle $C_R$ and at it's bound $W\geq 0$. It is followed that $W>0$ inside $C_R$
because of the extreme principle for harmonic functions.
Thus, inside $C_R$
\begin{equation}\label{29}
\frac{\omega}{2\pi}\int\limits_{D}\!\!\!\!\int
G(x,y,\xi,\eta)~d\xi ~
d\eta>\frac{\omega}{2\pi}\int\limits_{C_R}\!\!\!\!\int
G_R(x,y,\xi,\eta)~d\xi ~d\eta.
\end{equation}
From the representation \eqref{11} of function $\Psi$ and inequality \eqref{29}
inside $C_R$
\begin{equation}\label{30}
\Psi\leq U= C-\frac{\omega}{2\pi}\int\limits_{C_R}\!\!\!\!\int
G_R(x,y,\xi,\eta)~d\xi ~d\eta,\quad C=\max\varphi(s).
\end{equation}
Function $U$ in circle $C_R$ satisfies Poisson's equation $\Delta
U=\omega$ and on it's bound equals $C>0$.
So it gives us a possibility to write the function in explicit
form:
$$
U=C+\frac{\omega}{4}\,(r^2-R^2),\ r^2=x^2+y^2.
$$
If
\begin{equation}\label{31}
\omega>L=\frac{4C}{R^2}
\end{equation}
then $U\leq 0$ в $C_{R_1}$,\quad $R_1=\sqrt{R^2-4C/\omega}$. From described above and
\eqref{30} ($\Psi\leq U$) it's is followed that $\Psi\leq 0$ in $C_{R_1}.$

Thus, we got a simple condition for $\omega$ \eqref{31} to consider dual problem (2), (3) instead of hardly controllable condition \eqref{9}. An inequality
\eqref{29} is true for any circle $C_{R}\subset D$. But evaluation of volume $\omega$ increases in this case .

Let's consider M.A. Goldshtik's statement: if some set
$B^*$ realizes functional's extremum then at bound $\Psi=0.$
Described proof is true for boundary points that are inside domain $D$.
But there couldn't be excluded the possibility of nonzero cut of bounds
$B^*$ and $D$ which have points where $\varphi>0.$
In this case $B^*$ doesn't give minimum to a functional
(proves similarly to $B^*=D$). The function $\Psi$ is not a solution of problem (1), (2).
Let's assume that there is a maximum at $B^*$ and in all points of
$B^*$ $\varphi>0$. Then function $\Psi$ is a solution of problem (2), (3).
If in some point $\Psi<0$ then excluding from $B^*$
a sufficiently small neighborhood of a point where $\Psi<0$ we increase the functional.

So that means that if don't constrain additional conditions in variational statement
then extremums of the functional could not give solution of considered problems.
As function $\varphi\geq 0$ at bound $\Gamma$
then domain $B$ (where extremum of the functional is reached and gives the problem (1), (3) solution)
must not adjoin to points at the bound $\Gamma$ where $\varphi>0$, what is opposite to solution of problem (2), (3).

Let's assume that $\varphi>0$ everywhere at $\Gamma.$
Then domain of negativeness $B$ of the solution of problem (1), (3) should strictly locate inside $D$.
Domain of positiveness of problem solution (2), (3)  couldn't be the whole domain $D$
and one of it's connected boundary components must be the whole domain $D$.
On described above examples of model problems in this assumptions
with variational method was derived solution of problem (2), (3) and two nontrivial solutions of problem (1), (3).
It should be mentioned that solution of problem (1), (3)
was searched with variational method in simply connected domains.
In this approach it's not allowed to exclude the neighborhoods of the points what leads to increase of connectedness.
Thus, in variational statement of problem (1), (3)
it should be mentioned that functional is set in simply-connected domains which strictly lied in domain $D$.
So second nontrivial solution could be derived in the maximum case.
To find solution of a dual problem (2),~(3) using variational method the functional should be set in domains (not necessarily simply-connected), demanding one of the connected components of the bound coincide with the whole bound of the domain $D.$

Consider variational statement of problem (1), (3) for additional conditions that described above.
The functional \eqref{12} rewrite in the following form
\begin{equation}\label{32}
 I(B)=q+2\omega\int\limits_{B}\!\!\!\int\Psi_0dxdy-
 \frac{\omega^2}{2\pi}\int\limits_{B}\!\!\!\int dxdy \int\limits_{B}\!\!\!\int G(x,y,\xi,\eta)d\xi d\eta.
\end{equation}
$$
 q=\int\limits_{\Gamma}\Psi_0\frac{\partial\Psi_0}{\partial n}ds.
$$
Representation \eqref{11} of function $\Psi(x,y)$ is taken into account here. Because of Green's formula for harmonic functions
$$
 q=\int\limits_{\Gamma}\Psi_0\frac{\partial\Psi_0}{\partial n}ds=
 \int\limits_{D}\!\!\!\int(\nabla\Psi_0)^2dxdy\geq 0.
$$
Let's identify
$$
 \int\limits_{B}\!\!\!\!\int\Psi_0dxdy=A(B)>0,\quad
 \frac{1}{2\pi}\int\limits_{B}\!\!\!\int dxdy \int\limits_{B}\!\!\!\int G(x,y,\xi,\eta)d\xi d\eta=Q(B)>0.
$$
Then
\begin{equation}\label{33}
 I(B)=q+2\omega A(B)-\omega^2 Q(B)
\end{equation}
is a square trinomial in regard to $\omega$. When
\begin{equation}\label{34}
 \omega>L_1=\frac{A(1+\sqrt{1+\frac{Qq}{A^2}})}{Q}
\end{equation}
the functional \eqref{32} is negative in set $B$. The set $B_1$
could coincide with the whole domain $D$. From \eqref{33}
it is led that for any fixed $\omega$ exists $B_1\subset D$ where
\begin{equation}\label{35}
 I(B_1)>q.
\end{equation}
For trivial solution $\Psi_0(x,y)$ we have $I=q$.

Now evaluate $L_1$. Using inequality \eqref{29} we get
\begin{equation}\label{36}
I(D)\leq q+2С\pi R_1^2\omega C-\frac{\omega^2}{2\pi}\int\limits_{C_R}\!\!\!\!\!\int dxdy
\int\limits_{C_R}\!\!\!\!\!\int G_R(x,y,\xi,\eta)d\xi d\eta <0.
\end{equation}
where $R$ is a radius of the area-largest circle $C_R$ inscribed in domain $D$ and $R_1$ is a radius of the area-smallest circle inscribed in domain $D$.
After calculating integrals in \eqref{36} the following inequation is derived
\begin{equation}\label{37}
q+2\pi R_1^2 C \omega-\frac{\pi R^4}{8}\omega^2<0,
\end{equation}
which could be true when
\begin{equation}\label{38}
\omega>L_1=\frac{8R_1^2 C\left(1+\sqrt{1+\frac{qR^4}{8\pi R_1^4 C^2}}\,\right)}{R^4}.
\end{equation}

Now we consider extremums of the functional.
Let's suppose that extremum is realized in the whole domain $D$.
As was mentioned above minimum is impossible.
Let $\omega$ satisfy inequality \eqref{34} when
$B=D$, then the functional is negative.  Because of \eqref{33} for such $\omega$
$B$ could be selected so that $I(B)>0$.
Hence and from
$I(D)<0$ it is followed that $D$ doesn't realize the maximum.

Let's assume that $B$ realizes absolute maximum of the functional and $B\subset D$.
This set exists as in the opposite case the functional equals $q$,
and because of \eqref{33} domain $B$ could be selected that $I(B)>q$.
As Goldshtik showed in [4] function $\Psi$ vanishes at bound $B$.
Then $\Psi<0$ in $B$, because $\Delta\Psi=\omega>0$.
When a small neighborhood is excluded from $B$ the functional increases. However it can't be done,
because the functional is considered in simple-connected domains.

 Thus existence of nontrivial solution of problem (1), (3) is proved.
 Another nontrivial solution was derived by Goldshtik
 and it gives the minimal value to the functional.

Now we consider a functional at the bounds where one of the simply-connected components
coincides with the bound $\Gamma$ in domain $D$.
Demand of simple connectedness is not constrained.
Let's show that in the maximum case when inequality \eqref{31} is true,
there is a solution of dual problem (2), (3).
It could's be a maximum in the whole domain $D$
Because of inequality \eqref{31} in $D$ exists a circle $C_r$ where $\Psi<0$.
As there is no restriction for simple connectedness a small neighborhood could be excluded from $C_r$
in order to increase the functional.
If there is a minimum in $D$, as was mentioned before,
when inequality \eqref{31} is true function $\Psi$ doesn't give a solution of a dual problem.

Let's assume that a maximum of the functional is realized in $B^*$ and $B^*\neq D$.
One of the connected components of the bound $B^*$ coincides with $\Gamma$.
If in some point $M_0\in B^*$ function $\Psi<0$ then a small neighborhood of point  $M_0$ where $\Psi<0$ could be selected
so that exclusion of this point from $B^*$ according to \eqref{13} increases the functional.
The function $\Psi$ is harmonic in $D\setminus\overline{B}$ and vanishes at bound. So
 $\Psi\equiv 0$ in $D\setminus\overline{B}$.

Thus, it was established using introduced by M.A.~Goldshtik variational method
that when inequality \eqref{38} is true, the problem (1), (3) has second nontrivial solution.
Moreover, existence of the solution of dual problem (2), (3) is assigned using variational method.
It was proved in case of $\varphi(s)>0$.
General case
$\varphi(s)\geq 0$ doesn't substantially complicate the proof.

\medskip

\emph{ This research is financially supported by The Russian Foundation for Basic Research (grant 11-01-00283). }


\maketitEnd  

\end{document}